\documentclass[11pt]{article}
\usepackage{amsmath}
\usepackage{amssymb}
\usepackage{graphicx}
\usepackage{hyperref}

\renewcommand{\baselinestretch}{1.2}
  \renewcommand{\arraystretch}{1.1}
\begin{document}

 \title{An Elementary Proof for Ljunggren Equation}

   \author{Zhengjun Cao$^{1}$, \quad Lihua Liu$^{2,*}$}
  \footnotetext{ $^1$Department of Mathematics, Shanghai University, Shanghai, 200444,  China.  \\
   $^2$Department of Mathematics, Shanghai Maritime University, Shanghai, 201306,  China. \  $^*$\,\textsf{liulh@shmtu.edu.cn} }

 \date{}\maketitle
\begin{quotation}
 \textbf{Abstract}. We present an elementary proof for Ljunggren equation.


  \end{quotation}

  \section{Introduction}
In 1942,  Ljunggren \cite{L42} presented a sophisticated proof for that the only solutions of  equation
$$1+x^2=2y^4  \eqno(1)$$  in positive integers are $(1, 1)$ and $(239, 13)$. Mordell \cite{G04} asked  if it was possible to find a simple or elementary proof.
In 1991, Steiner and Tzanakis \cite{ST91} gave a simpler proof by reducing the problem to a Thue equation and then solving it by using a deep result of Mignotte and Waldschmidt \cite{MW88} on linear forms in logarithms and continued fractions. In this not, we present an elementary proof for the problem.

\section{Preliminaries}

\textbf{Lemma 1}. \emph{For any two adjacent numbers $a, b$ in sequence $\{a_i\}_{i=0}^{\infty}$ where $a_{0}=0, a_1=1, a_{n+2}=2a_{n+1}+a_{n}$, $a^2+b^2$ is just in the sequence.}

\emph{Proof.} The  sequence is 
$$0,  1, 2, 5, 12, 29, 70, 169, 408, 985, 2378, 5741, 13860, 33461, \cdots \eqno(I) $$
For any two adjacent numbers $a, b$ in $(I)$, we obtain the  following sub-sequence
$$a, b, a+2b, 2a+5b, 5a+12b, 12a+29b, 29a+70b, 70a+169b, \cdots   $$
 Apparently, the coefficients of $b$ consist of the original sequence $(I)$. Therefore,
   $a^2+b^2$ is just in the sequence. \hfill $\Box$

   For example,
$0^2+1^2=1, 1^2+2^2= 5, 2^2+5^2=29, 5^2+12^2=169, $
$12^2+ 29^2=985, 29^2+70^2=5741, 70^2+169^2=33461$. 
Immediately, we have the following corollary.

\textbf{Corollary 1}. \emph{Each number in the following sub-sequence,
$$1,\ 5,\ 29,\ 169,\ 985,\ 5741,\ 33461,\ \cdots  \eqno(II)$$
can be represented as the sum of two squares of two adjacent numbers in sequence $(I)$.}

\emph{Proof.} In the sequence 
$$0, 1, \cdots, a, b, 2b+a, \cdots, \underline{a\cdot a+ b\cdot b}, b\cdot a+(2b+a)b, \underline{(2b+a)a+(5b+2a)b}, \cdots   $$
we have $b^2+(2b+a)^2=a^2+4ab+5b^2=(2b+a)a+(5b+2a)b$.

\textbf{Lemma 2}. \emph{For any two adjacent numbers $a, b$ in sequence $(I)$, if $b\neq 1$ or $12$, then
$a^2+b^2$ is not a square.}

\emph{Proof.}  In sequence $(I)$, \emph{every segment of length 5}  can be represented as follows
$$c, d, c+2d, 2c+5d, 5c+12d $$
If $b>12$ then $a, b$ can be represented as $2c+5d, 5c+12d$ ($c\neq 0$), respectively.  Hence,
$$a^2+b^2=(2c+5d)^2+(5c+12d)^2=169d^2+140cd+29c^2 $$

If there exists a positive integer $k $ such that
$$ 169d^2+140cd+29c^2=k^2 \eqno(2) $$
then  we have
$$ 169\left(\frac{d}{c}\right)^2+140\left(\frac{d}{c}\right)+29-\left(\frac{k}{c}\right)^2=0 \eqno(3)$$
Its  discriminant is
$$\triangle=4\left(\frac{169 k^2-c^2}{c^2}\right) $$

Suppose that
$$ 169 k^2-c^2= \beta^2 \eqno(4)$$
for some positive  integer $\beta$. Hence, $(c, \beta, 13k)$ is a solution of equation $X^2+Y^2=Z^2$.

  If $13\nmid c$, then $\beta/c =5/12$ or $12/5$. By Eq.(3), we have
$$ \frac{d}{c}=\frac{-140+\frac{2\beta}{c}}{2\times 169 }<0 \eqno(5)$$
It is a contradiction.

  If $13 |c$, then by Eq.(4), $13 |\beta$,  and  by Eq.(2),  $13|k.$  Eventually,  we have
 $$ 169 \hat{k}^2-\hat{c}^2= \hat{\beta}^2 \eqno(4')$$
 for some positive integers $\hat{k}, \hat{c}, \hat{\beta}$, where $(\hat{\beta}, \hat{c})=1$ and  $13\nmid\hat{c}$. Thus, $\beta/c=\hat{\beta}/\hat{c} =5/12$ or $12/5$. By Eq.(5), it leads to a contradiction, too.
  \hfill $\Box$

\section{An elementary proof}

\textbf{Theorem 1}. \emph{The only solutions of
$1+x^2=2y^4 $  in positive integers are $(1, 1)$ and $(239, 13)$.}

\emph{Proof.}  It suffices to consider the following Pell's equation
$$x^2-2(y^2)^2=-1 $$
Since $|\frac{x}{y^2}-\sqrt{2}|=|\frac{1}{y^4}\times \frac{1}{\frac{x}{y^2}+\sqrt{2}} |<\frac{1}{2y^4} $,  $\frac{x}{y^2}$ is a convergent of the regular continued fraction for $\sqrt{2}$. Note that  $\sqrt 2= [1; 2, 2, 2, \cdots] $.  Its convergents are
$$ 1,\ \frac{3}{2}, \ \frac{7}{5},\ \frac{17}{12},\ \frac{41}{29},\ \frac{99}{70},\ \frac{239}{169}, \ \frac{577}{408},\ \frac{1393}{985}, \ \frac{3363}{2378}, \frac{8119}{5741}, \frac{19601}{13860}, \frac{47321}{33461}, \cdots  $$
The solutions of equation $X^2-2Y^2=-1 $
are the convergents
$$ 1,\  \frac{7}{5},\  \frac{41}{29},\  \frac{239}{169}, \  \frac{1393}{985}, \  \frac{8119}{5741},  \frac{47321}{33461},  \cdots $$
Hence,  the original problem is reduced  to finding squares in the sequence
$$1,\ 5,\ 29,\ 169,\ 985,\ 5741,\ 33461,\ \cdots  \eqno(II)$$

By Corollary 1, each number in $(II)$ can be represented as $a^2+b^2$ for some two adjacent numbers $a, b$ in $(I)$.
Furthermore, by Lemma 2, there are only two squares $1, 169$ in sequence $(II)$. Thus, $y=1,$ or $13$.
 \hfill $\Box$

\section{Conclusion}
We present an elementary proof for Ljunggren equation which solves Mordell's open problem.

%\textbf{Acknowledgments}. We thank the National Natural Science Foundation of China (Project
%61303200, 61411146001).

\end{document}